# Explicit formulae for spectral norms of circulant-type matrices with some given entries

J.W.Zhou*, Z.L.Jiang

Department of Mathematics, Linyi University, China

*Corresponding author: jwzhou@yahoo.com

**Abstract**

In this paper we investigate the spectral norm for circulant matrices, whose entries are modified Fibonacci numbers and Lucas numbers. We obtain the identity estimations for the spectral norms. Some numerical test results are listed to verify the results using those approaches.

**Introduction**

The Fibonacci and Lucas sequences $\{F_n\}$ and $\{L_n\}$ are defined by the recurrence relations

$$F_0 = 0, \quad F_1 = 1, \quad \cdots, \quad F_n = F_{n-1} + F_{n-2}, \quad \cdots,$$

and

$$L_0 = 2, \quad L_1 = 1, \quad \cdots, \quad L_n = L_{n-1} + L_{n-2}, \quad \cdots.$$

If we start from $n = 0$, then the Fibonacci and Lucas sequences are given by

$$
\begin{array}{ccccccccc}
n & 0 & 1 & 2 & 3 & 4 & 5 & 6 & 7 & \ldots \\
F_n & 0 & 1 & 1 & 2 & 3 & 5 & 8 & 13 & \ldots \\
L_n & 2 & 1 & 3 & 4 & 7 & 11 & 18 & 29 & \ldots.
\end{array}
$$

In *(Kocer, Tuglu and Stakhov, 2009)*, Their Binet forms are given by

$$F_n = \frac{1}{\sqrt{5}}\left[\left(\frac{1+\sqrt{5}}{2}\right)^n - \cos(\pi n)\left(\frac{1+\sqrt{5}}{2}\right)^{-n}\right],$$

$$L_n = \left(\frac{1+\sqrt{5}}{2}\right)^n + \cos(\pi n)\left(\frac{1+\sqrt{5}}{2}\right)^{-n}.$$

There is a large class of identities in the construction as the following two tables. The following identities for the Fibonacci numbers are well known in previous work by the authors (Akbulak and Bozkurt, 2008; Benoumhani, 2003; Melham, 1999; Ipek, 2011):

Table 1. The Fibonacci number identities of type

| $No.$ | $a_i$ | $s_n$ | $No.$ | $a_i$ | $s_n$ |
|---|---|---|---|---|---|
| 1 | $F_i$ | $F_{n+1} - 1$ | 5 | $F_i F_{n-1-i}$ | $\frac{nF_{n+1}+(n+2)F_{n+1}-2}{5}$ |
| 2 | $F_{2i+1}$ | $F_{2n}$ | 6 | $\sum_{k=0}^{i} F_k^2$ | $F_n^2 + \frac{(-1)^n - 1}{2}$ |
| 3 | $F_{2i}$ | $F_{2n-1} - 1$ | 7 | $F_i F_{i+1}$ | $F_n^2 + \frac{(-1)^n - 1}{2}$ |
| 4 | $iF_i$ | $(n-1)F_{n+1} - F_{n+2} + 2$ | | | |

Applying $L_{n+1} = F_n + F_{n+2}$ and the identities of Fibonacci numbers in above table, we deduce the sum formulates of Lucas numbers as following:

Table 2. The Lucas number identities of type $\sum_{i=0}^{n-1} a_i = s_n$

| No. | $a_i$ | $s_n$ | No. | $a_i$ | $s_n$ |
|---|---|---|---|---|---|
| 1 | $L_i$ | $L_{n+1} - 1$ | 2 | $L_{2i+1}$ | $2F_{2n+1} - F_{2n} - 2$ |
| 3 | $L_{2(i+1)}$ | $F_{2n+1} + 2F_{2n} - 1$ | 4 | $F_i L_i$ | $2F_{n-1}^2 + F_{n-1}F_n + (-1)^{n-1} - 1$ |
| 5 | $L_i F_{n-1-i}$ | $nF_n$ | 6 | $F_i + L_i$ | $4F_n + 2F_{n-1} - 2$ |
| 7 | $\sum_{k=0}^{i} L_k^2$ | $L_n^2 + 2n + (-1)^{n-1}\frac{5}{2} - \frac{3}{2}$ | 8 | $L_i L_{i+1}$ | $L_n^2 + (-1)^{n-1}\frac{5}{2} - \frac{3}{2}$ |
| 9 | $iF_i$ | $(3n-5)F_n + (n-2)F_{n-1} - 2F_{n+1} + 4$ | | | |

**Definition** 1. *[Stallings and Boullion,1972]* A circulant matrix is an $n \times n$ complex matrix with the following form:

$$A = \begin{pmatrix} a_0 & a_1 & \cdots & a_{n-1} \\ a_{n-1} & a_0 & \cdots & a_{n-2} \\ a_{n-2} & a_{n-1} & \cdots & a_{n-3} \\ \vdots & \vdots & \ddots & \vdots \\ a_1 & a_2 & \cdots & a_0 \end{pmatrix}_{n \times n}. \quad (1)$$

The first row of $A$ is $(a_0, a_1, \cdots, a_{n-1})$, its $(j+1)$-th row is obtained by giving its $j$-th row a right circular shift by one positions.

**Definition** 2. *[Akbulak and Bozkurt,2008]* The spectral norm $\|\cdot\|_2$ of a matrix $A$ with complex entries is the square root of the largest eigen-value of the positive semi-definite matrix $A^*A$:

$$\|A\|_2 = \sqrt{\lambda_{max}(A^*A)},$$

where $A^*$ denotes the conjugate transpose of $A$. Therefore if $A$ is an $n \times n$ real symmetric matrix or $A$ is a normal matrix, then

$$\|A\|_2 = \max_{1 \leq i \leq n} |\lambda_i|,$$

where $\lambda_1, \lambda_2, \cdots, \lambda_n$ are the eigen-values of $A$.

**The identity estimations for spectral norms**

We give the main theorems of this paper in the following parts.

**Theorem 1.** Set $B_1$ as the matrix in Eq.(1), the first row of $B_1$ is $(0F_0, 1F_1, \ldots, (n-1)F_{n-1})$, then we have the following identity:

$$\|B_1\|_2 = (n-1)F_{n+1} - F_{n+2} + 2.$$

**Proof.** As we all known, $B_1$ is normal. Combining with **Definition 2**, the spectral radius of $B_1$ is equal to its spectral norm. Furthermore, applying the irreducible and entry-wise nonnegative properties, we claim that $\|B_1\|_2$ is equal to its Perron value.

We set an $n$-dimensional column vector $v = (1, 1, \cdots, 1)^T$, then

$$B_1 v = \left(\sum_{i=0}^{n-1} iF_i\right) v.$$

Obviously, $\sum_{i=0}^{n-1} iF_i$ is an eigen-value of $B_1$ associated with the positive eigen-vector $v$, which is the Perron value of $B_1$. Employing the identities of Fibonacci numbers in **Table 1**, we have

$$\|B_1\|_2 = (n-1)F_{n+1} - F_{n+2} + 2.$$

This completes the proof.

**Theorem 2.** Let $B_2$ with the form as Eq.(1), the first row of $B_2$ is $(F_0 F_1, F_1 F_2, \ldots, F_{n-1} F_n)$, then we have the following identity:

$$\|B_2\|_2 = F_n^2 + \frac{(-1)^n - 1}{2}.$$

**Proof.** Using the same techniques of **Theorem 1**, the normal matrix $B_2$ is irreducible and entry-wise nonnegative, so, the spectral norm of $B_2$ is the same as its Perron value. Let

$$v^T = \underbrace{(1, 1, \cdots, 1)}_{n},$$

then $B_2 v = \left(\sum_{i=0}^{n-1} F_i F_{i+1}\right)v$. As $\sum_{i=0}^{n-1} F_i F_{i+1}$ is an eigen-value of $B_2$ associated with the positive eigen-vector $v$, which is equal to Perron value of $B_2$. Combining with the identities of Fibonacci numbers in **Table 1**, we obtain

$$\|B_2\|_2 = F_{n-1} F_n,$$

which completes the proof.

**Lemma 1.**[*Horn and Johnson*, 1986] Let $A$ is a nonnegative matrix, if the column sums of $A$ are constant, then

$$\rho(A) = \|A\|_1,$$

where $\rho(A) = \max\{|\lambda| : \lambda \text{ is an eigenvalue of matrix } A\}$, and $\|\cdot\|_1$ denotes the maximum column sum matrix norm.

**Theorem 3.** Let $B_3$ with the form as Eq.(1), and the first row of $B_3$ is $(F_1, F_3, \ldots, F_{2n-1})$, then we have the following identity:

$$\|B_3\|_2 = F_{2n}.$$

**Proof.** Obviously, $B_3$ is normal. By **Definition 2**, we declare that the spectral radius of $B_3$ is equal to $\rho(B_3)$, $i.e.$, $\|B_3\|_2 = \rho(B_3)$. Furthermore, applying entry-wise nonnegative properties and column sums of $B_3$ are constant. By **Lemma 1**, we obtain

$$\rho(B_3) = \|B_3\|_1.$$

Employing the identities of Fibonacci numbers in **Table 1**, we have

$$\|B_3\|_2 = F_{2n}.$$

Furthermore, we give some theorems without proofs, which can be proved with the same approaches as the above theorems.

**Theorem 4.** Let some circulant matrices as Eq.(1), we obtain
(1) if the first row of $B_4$ is $(F_0, F_2, \ldots, F_{2n-2})$, then

$$\|B_4\|_2 = F_{2n-1} - 1.$$

(2) if the first row of $B_5$ is $(F_0 F_{n-1}, F_1 F_{n-2}, \ldots, F_{n-1} F_0)$, then

$$\|B_5\|_2 = \frac{nF_{n+1} + (n+2)F_{n+1} - 2}{5}.$$

(3) if the first row of $B_6$ is $(0L_0, 1L_1, \ldots, (n-1)L_{n-1})$, then

$$\|B_6\|_2 = (3n-5)F_n + (n-2)F_{n-1} - 2F_{n+1} + 4.$$

(4) if the first row of $B_7$ is $\left(L_0^2, L_1^2, \ldots, L_{n-1}^2\right)$, then

$$\|B_7\|_2 = L_{n-1} L_n + 2.$$

(5) if the first row of $B_8$ is $(L_1, L_3, \ldots, L_{2n-1})$, then

$$\|B_8\|_2 = 2F_{2n+1} - F_{2n} - 2.$$

(6) if the first row of $B_9$ is $(L_0, L_2, \ldots, L_{2n-2})$, then

$$\|B_9\|_2 = F_{2n+1} + 2F_{2n} - 1.$$

(7) if the first row of $B_{10}$ is $(L_0 L_1, L_1 L_2, \ldots, L_{n-1} L_n)$, then

$$\|B_{10}\|_2 = \begin{cases} L_n^2 - 4 & n \text{ even,} \\ L_n^2 + 1 & n \text{ odd.} \end{cases}.$$

(8) if the first row of $B_{11}$ is $\left(L_0^2, \sum_{k=0}^{1} L_k^2, \ldots, \sum_{k=0}^{n-1} L_k^2\right)$, then

$$\|B_{11}\|_2 = \begin{cases} L_n^2 + 2n - 4 & n \text{ even,} \\ L_n^2 + 2n + 1 & n \text{ odd.} \end{cases}.$$

(9) if the first row of $B_{12}$ is $(F_0 L_0, F_1 L_1, \ldots, F_{n-1} L_{n-1})$, then

$$\|B_{12}\|_2 = 2F_{n-1}^2 + F_{n-1} F_n + (-1)^{n-1} - 1.$$

(10) if the first row of $B_{13}$ is $(F_0 + L_0, F_1 + L_1, \ldots, F_{n-1} + L_{n-1})$, then

$$\|B_{13}\|_2 = 4F_n + 2F_{n-1} - 2.$$

(11) if the first row of $B_{14}$ is $(L_0 F_{n-1}, L_1 F_{n-2}, \ldots, L_{n-1} F_0)$, then

$$\|B_{14}\|_2 = nF_n.$$

**Numerical examples**

**Example.** In this example, we give the numerical results for $B_1 - B_{14}$ in following Table 3.

Table 3. Spectral norms of $B_1 - B_{14}$

| $n$ | $B_1$ | $B_2$ | $B_3$ | $B_4$ | $B_5$ | $B_6$ | $B_7$ | $B_8$ | $B_9$ | $B_{10}$ | $B_{11}$ | $B_{12}$ | $B_{13}$ | $B_{14}$ |
|---|---|---|---|---|---|---|---|---|---|---|---|---|---|---|
| 3 | 3 | 3 | 4 | 8 | 1 | 7 | 14 | 12 | 16 | 17 | 23 | 4 | 8 | 3 |
| 5 | 21 | 24 | 33 | 55 | 5 | 47 | 79 | 77 | 121 | 122 | 132 | 33 | 24 | 15 |
| 6 | 46 | 64 | 88 | 144 | 10 | 102 | 200 | 200 | 320 | 320 | 332 | 88 | 40 | 30 |
| 7 | 94 | 168 | 232 | 377 | 20 | 210 | 524 | 522 | 841 | 842 | 856 | 232 | 66 | 56 |
| 8 | 185 | 441 | 609 | 987 | 38 | 413 | 1365 | 1365 | 2205 | 2205 | 2221 | 609 | 108 | 104 |

The above results demonstrate that the identities of spectral norms of $B_1 - B_{14}$ hold.

**Conclusions**

This paper had discussed the identity estimates of spectral norms for some circulant matrices, which are listed by explicit formulations. It is an interesting problem to investigate properties for circulant matrices with certain modified Fibonacci numbers and Lucas numbers, including norms, determinants, inverses, and so on.

**Acknowledgements**

The research was supported by National Natural Science Foundation of China (Grant No.11201212), Promotive Research Fund for Excellent Young and Middle-aged Scientists of Shandong Province (Grant No. BS2012DX004), and the Special Funds for Doctoral Authorities of Linyi University.

**References**

Akbulak M., Bozkurt D. (2008), On the norms of Toeplitz matrices involving Fibonacci and Lucas numbers. *Hacet. J. Math. Stat., 37*, pp. 89–95.
Benoumhani M. (2003), A Sequence of Binomial Coefficients Related to Lucas and Fibonacci Numbers. *Journal of Integer Sequences, Portugaliae Mathematica, 6*, pp. 1–10.
Bose A., Hazra R. S. and Saha K. (2011), Poisson convergence of eigen-values of circulant type matrices. *Extremes, 14*, pp. 365–392.
Bose A., Hazra R. S. and Saha K. (2011), Spectral norm of circulant-type matrices. *J. Theor. Probab., 24*, pp. 479–516.
Bose A., Guha S., *et al*.(2011), Circulant type matrices with heavy tailed entries. *Statistics and Probability Letters, 81*, pp. 1706–1716.
Chou W. S., Du B. S., *et al*.(2008), A note on circulant transition matrices in Markov chains. *Linear Algebra Appl., 429*, pp. 1699–1704.
Erbas C., Tanik M. M. (1995), Generating solutions to the N-Queens problems using 2-circulants. *Math. Mag., 68*, pp. 343–356.
Kocer E. G., Tuglu N. and Stakhov A. (2009), On the $m$-extension of the Fibonacci and Lucas $p$-numbers. *Chaos Soliton. Fract., 40*, pp. 1890–1906.
Horn R.A., Johnson C.R.(1986), Matrix Analysis. *Cambridge: Cambridge University Press.*
Melham R. (1999), Sums involving Fibonacci and Pell numbers. *Portugaliae Mathematica, 56 (3)*, pp. 1–9.

Ngondiep E., Stefano S. C. and Sesana D. (2010), Spectral features and asymptotic properties for $g$-circulants and $g$-Toeplitz sequences. *SIAM J. Matrix Anal. Appl., 31*, pp. 1663–1687.

Stallings W. T. and Boullion T. L. (1972), The pseudoinverse of an $r$–circulant matrix, *Proc. AMS.,34(2)*, pp. 385–388.

Ipek A.(2011), On the spectral norms of circulant matrics with classical Fibonacci and Lucas numbers entries. *Appl. Math. Comput., 217*, pp. 6011–6012.

Woods N. A., Galatsanos N. P. and Katsaggelos A. K. (2006), Stochastic methods for joint registration, restoration, and interpolation of multiple undersampled images. *Linear Algebra Appl., 15*, pp. 201–213.

Wu Y. K., Jia R. Z. and Li Q. (2002), $g$-circulant solutions to the $(0,1)$ matrix equation $A^m = J_n^\star$. *Linear Algebra Appl., 345*, pp. 195–224.